\documentclass[12pt]{article}
\usepackage{amsfonts, amsmath, amstext, amssymb}

\usepackage[dvips]{graphicx}
\usepackage{latexsym}

\newtheorem{thm}{Theorem}[section]
\newtheorem{prop}[thm]{Proposition}

\newtheorem{lemma}[thm]{Lemma}

\newtheorem{cor}[thm]{Corollary}

\newtheorem{definitiontemp}[thm]{Definition}
\newenvironment{defn}{\begin{definitiontemp}
\normalfont}{\end{definitiontemp}}

\def\bec{\begin{cor}}
\def\enc{\end{cor}}

\def\ep{\varepsilon}

\def\bet{\begin{thm}}
\def\ent{\end{thm}}

\def\becor{\begin{cor}}
\def\encor{\end{cor}}

\def\bel{\begin{lem}}
\def\enl{\end{lem}}

\def\bedef{\begin{defn}}
\def\endef{\end{defn}}

\def\bep{\begin{prop}}
\def\enp{\end{prop}}

\newenvironment{pf}{\begin{trivlist}\item[\hskip\labelsep
{\it Proof.}]}{\end{trivlist}}

\newcommand{\Fin}{\textbf{Fin}}

\newcommand{\Inf}{\textbf{Inf}}

\newcommand{\set}[2]{\ensuremath{ \{ #1 : #2 \} }}

\renewcommand{\deg}[1]{\ensuremath{\text{deg}(#1)}}

\newcommand{\N}{\mathbb{N}}
\newcommand{\Z}{\mathbb{Z}}

\newcommand{\Q}{\mathbb{Q}}
\newcommand{\R}{\mathcal{R}}
\newcommand{\C}{\mathcal{C}}
\newcommand{\A}{\mathcal{A}}
\newcommand{\B}{\mathcal{B}}

\renewcommand{\P}{\mathbb{P}}
\newcommand{\U}{\mathcal{U}}

\newcommand{\Ubar}{\overline{U}}

\newcommand{\Wbar}{\overline{W}}

\newcommand{\Uvec}{\vec{U}}
\newcommand{\Vvec}{\vec{V}}

\newcommand{\xvec}{\vec{x}}
\newcommand{\Tvec}{\vec{T}}
\newcommand{\Xvec}{\vec{X}}

\newcommand{\Yvec}{\vec{Y}}

\newcommand{\Zvec}{\vec{Z}}

\newcommand{\gp}{{\mathfrak p}}

\newcommand{\GL}[1]{\textbf{GL}_{#1}}

\newcommand{\la}{\langle}
\newcommand{\ra}{\rangle}

\newcommand{\Int}[1]{\text{Int}(#1)}

\newcommand{\SubQ}{\textbf{Sub}(\Q)}
\newcommand{\HTP}[1]{HTP(#1)}
\newcommand{\HTPQ}{HTP(\Q)}
\newcommand{\HTPZ}{HTP(\Z)}

\def\diverges{\!\uparrow}
\def\converges{\!\downarrow}
\newcommand{\at}{\char'100}

\newcommand{\cO}{{\mathcal O}}
\newcommand{\Leg}[2]{\left(\begin{smallmatrix} \underline{#1} \\ #2 \end{smallmatrix}\right)} 
\newcommand{\ov}[1]{\overline{#1}}

\newcommand{\qed}{\hbox to 0pt{}\nobreak\hfill\rule{2mm}{2mm}}

\newcommand{\bfc}{\boldsymbol{c}}

\def\bfz{\boldsymbol{0}}
\def\s01{\ensuremath{\Sigma^0_1}}
\def\d02{\ensuremath{\Delta^0_2}}
\def\phi{\varphi}

\def\res{\!\!\upharpoonright\!}

\newcommand{\comment}[1]{}
\renewcommand{\omega}{\mathbb{N}}

\begin{document}

\title{The Hilbert's-Tenth-Problem Operator}
\author{Kenneth Kramer \& Russell Miller
\thanks{Research of the first author was partially supported by a PSC-CUNY award, cycle 48,
jointly funded by the Professional Staff Congress and C.U.N.Y.  The second 
author was partially supported by Grant \# DMS -- 1362206
from the National Science Foundation, and by several grants
from the PSC-CUNY Research Award Program.
}
}
\maketitle

\begin{abstract}
For a ring $R$, Hilbert's Tenth Problem $HTP(R)$ is the set of polynomial
equations over $R$, in several variables, with solutions in $R$.
We view $HTP$ as an operator, mapping each set $W$ of prime numbers
to $HTP(\Z[W^{-1}])$, which is naturally viewed as a set of polynomials
in $\Z[X_1,X_2,\ldots]$.  For $W=\emptyset$, it is a famous result
of Matiyasevich, Davis, Putnam, and Robinson that the jump
$\emptyset'$ is Turing-equivalent to $HTP(\Z)$.  More generally,
$HTP(\Z[W^{-1}])$ is always Turing-reducible to $W'$,
but not necessarily equivalent.  We show here 
that the situation with $W=\emptyset$ is anomalous:
for almost all $W$, the jump $W'$ is not diophantine in $HTP(\Z[W^{-1}])$.
We also show that
the $HTP$ operator does not preserve Turing equivalence:
even for complementary sets $U$ and $\Ubar$,
$HTP(\Z[U^{-1}])$ and $HTP(\Z[\Ubar^{-1}])$ can differ by a
full jump.  Strikingly, reversals are also possible, with
$V<_T W$ but $\HTP{\Z[W^{-1}]} <_T \HTP{\Z[V^{-1}]}$.
\end{abstract}

\section{Introduction}
\label{sec:intro}

The original version of Hilbert's Tenth Problem demanded an algorithm
deciding which polynomial equations from $\Z[X_0,X_1,\ldots]$ have solutions
in integers.  In 1970, Matiyasevic \cite{M70} completed work by Davis, Putnam
and Robinson \cite{DPR61}, showing that no such algorithm exists.
In particular, these authors showed that there exists a $1$-reduction from the
Halting Problem $\emptyset'$ to the set of such equations with solutions,
by proving the existence of a single polynomial $h\in\Z[X,\Yvec]$
such that, for each $n$ from the set $\omega$
of nonnegative integers, the polynomial $h(n,\Yvec)=0$ has a solution
in $\Z$ if and only if $n$ lies in $\emptyset'$.
Since the membership in the Halting Problem was known
to be undecidable, it followed that Hilbert's Tenth Problem was also undecidable.

One naturally generalizes this problem to all rings $R$, defining
Hilbert's Tenth Problem for $R$ to be the set
$$ \HTP{R} = \set{f\in R[\Xvec]}{(\exists r_1,\ldots,r_n\in R^{<\omega})~f(r_1,\ldots,r_n)=0}.$$
Here we will examine this problem for one particular class:  the subrings $R$
of the field $\Q$ of rational numbers.  Notice that in this situation,
deciding membership in $\HTP{R}$ reduces to the question of deciding
this membership just for polynomials from $\Z[\Xvec]$, since one
readily eliminates denominators from the coefficients of a polynomial in $R[\Xvec]$.
So, for us, $\HTP{R}$ will always be a subset of $\Z[X_1,X_2,\ldots]$.  In turn,
sets of polynomials, such as $\HTP{R}$, will be viewed as subsets of $\omega$,
using a fixed computable bijection from $\omega$ onto $\Z[\Xvec] = \Z[X_0,X_1,\ldots]$.

Subrings $R$ of $\Q$ correspond bijectively to subsets $W$ of the set $\P$
of all primes, via the map $W\mapsto\Z[\frac1p:p\in W]$.  We write $R_W$
for the subring $\Z[\frac1p:p\in W]$, which (as a subset of $\Q$) is Turing-equivalent
to $W$ itself (as a subset of $\P$). The \emph{$HTP$ operator} is the map
sending each $W\subseteq\P$ to $\HTP{R_W}$.  This operator, and its
relation to Turing reducibility $\leq_T$, are the focus of our work here.  
Recall that $A\leq_T B$ intuitively means that, if one knew which
numbers lie in $B$, one could decide which numbers lie in $A$.
($A<_T B$ just means that $A\leq_T B$ but $B\not\leq_T A$.)  One would
expect that, for more complicated rings $R$, $\HTP{R}$ would be more
difficult to compute than for simpler rings.  We will confound this expectation,
employing properties of the polynomials $X^2+qY^2-1$ to produce
subrings $R$ and $S$ of $\Q$ such that $R<_T S$ -- so it is strictly
easier to decide which rationals lie in $R$ than which lie in $S$ --
yet $\HTP{S}<_T\HTP{R}$ -- i.e., it is strictly easier to decide which
polynomials have solutions in $S$ than to decide which have solutions in $R$.

Of course, we have a computable bijection between subsets of $\omega$
and subsets of $\P$, using the function mapping $n\in\omega$ to the $n$-th prime $p_n$,
starting with $p_0=2$.  Since this bijection preserves Turing degrees,
Turing reductions can use sets $W$ in $2^{\P}$ as their oracles
either by converting them to subsets of $\omega$ or by viewing them
as subsets of $\omega$ already (which happen to be subsets of $\P$).
The $HTP$ operator defined above specifically uses a subset of $\P$ as its input,
so we will generally stick to subsets of $\P$ in our notation in this article.

Occasionally we will consider the jump $\HTP{R}'$ of a set $\HTP{R}$,
using the conversion just described.  Recall that the \emph{jump} $W'$
of a set $W\subseteq\omega$ is essentially the Halting Problem,
relativized to $W$:
$$ W' = \set{e\in\omega}{\Phi_e^W(e)\text{~halts}},$$
where $\Phi_e$ is the $e$-th \emph{oracle Turing machine},
and $\Phi_e^W$ denotes the partial function computed by this machine
when it runs with $W$ as its oracle.  (Details about jumps, Turing reducibility,
and oracle Turing computation
may be found in many standard sources, such as \cite[Chap.\ III]{S87}.)
\comment{
It is well known that, whenever $A\equiv_T B$, we must also have $A'\equiv_T B'$.
Therefore, the jump $W'$ of a subset $W\subseteq\P$ has the same Turing degree
as the jump of its preimage $\set{n\in\omega}{p_n\in W}\subseteq\omega$ under
the bijection $n\mapsto p_n$.  Therefore, we will usually apply
the jump operator without worrying whether the given set is a subset of $\P$
or just of $\omega$.
}

One normally views subsets of $\P$ as paths through the tree $2^{<\P}$,
a complete binary tree whose nodes are the functions from initial segments
of the set $\P$ into the set $\{ 0,1\}$.  This allows us to introduce a topology on the space
$2^{\P}$ of paths through $2^{<\P}$, and thus on the class $\SubQ$ of all subrings of $\Q$.
Each basic open set $\U_\sigma$ in this topology is described by a node $\sigma$ on the tree:
$\U_\sigma = \set{W\subseteq\P}{\sigma\subset W}$,
where $\sigma\subset W$ denotes that when $W$ is viewed as a function
from $\P$ into the set $2=\{ 0,1\}$ (i.e., as an infinite binary sequence),
$\sigma$ is an initial segment of that sequence.
Also, we put a natural measure $\mu$ on $\SubQ$:
just transfer to $\SubQ$ the obvious Lebesgue measure on the power set $2^{\P}$ of $\P$.
Thus, if we imagine choosing a subring $R$ by flipping a fair coin (independently for
each prime $p$) to decide whether $\frac1p\in R$, the \emph{measure} of a subclass
$\mathcal S$ of $\SubQ$ is the probability that the resulting subring will lie in $\mathcal S$.

It is also natural, and in certain respects more productive, to consider Baire category
theory on the space $\SubQ$, as an alternative to measure theory.
For background regarding Baire category theory
on subrings of $\Q$, we refer the reader to \cite{MCiE16},
while parallel discussion of measure theory occurs in \cite{MIMS}.
Due to the common subject matter of those articles and this one, there is a
substantial overlap between the introductions and background sections
of the three papers, which we trust the reader to forgive.  Naturally,
we have also made every effort to maintain the same notation across
both papers.

\section{Subrings of the Rationals}
\label{sec:subrings}

Now we recall certain specific results about subrings of $\Q$.
For all $W\subseteq \P$, writing $R_W$ for $\Z[W^{-1}]$ as before,
we have the Turing reductions
$$ W\oplus\HTPQ \leq_T \HTP{R_W} \leq_T W'.$$
Indeed, each of these two Turing reductions is a \emph{$1$-reduction};
details appear in \cite[\S 2.2]{MIMS}.
Recall that the \emph{semilocal} subrings of $\Q$ are precisely those of the form
$R_W$ where the set $W$ is cofinite in $\P$, containing all but finitely many primes.
It will be important for us to know that whenever
$R$ is a semilocal subring of $\Q$, we have $\HTP{R}\leq_1\HTPQ$.
Indeed, both the Turing reduction and the $1$-reduction are uniform in the complement.
This result, which follows from Corollary \ref{cor:semilocal} below,
began with work of Julia Robinson in \cite{R49}.  A proof by
Eisentr\"ager, Park, Shlapentokh, and the author appears in \cite{EMPS15},
based in turn on work by Koenigsmann in \cite{K15}.
\begin{prop}[see Proposition 5.4 in \cite{EMPS15}]
\label{prop:semilocal}
For every prime $p$, there is
a polynomial $g_{p}(Z,X_1,X_2,X_3)$ such that for all rationals $q$, we have
$$ q\in R_{(\P-\{p\})} \iff g_{p}(q,\Xvec)\in\HTPQ.$$
Moreover, $g_p$ may be computed uniformly in $p$.
\end{prop}
\begin{cor}
\label{cor:semilocal}
For each finite subset $A_0\subseteq\P$, a polynomial
$f(Z_0,\ldots,Z_n)$ has a solution in $R_{(\P-A_0)}$ if and only if
$$(f(\Zvec))^2 + \sum_{p\in A_0, j\leq n} (g_{p}(Z_j,X_{1,j,p},X_{2,j,p},X_{3,j,p}))^2$$
has a solution in $\Q$.
\end{cor}

\comment{
There exists a computable function $G$ such that
for every $n$, every finite set $A_0=\{ p_1,\ldots,p_n\}\subset\P$ and every $f\in\Z[\Xvec]$,
$$ f\in\HTP{R_{(\P-A_0)}} \iff G(f,\la p_1,\ldots,p_n\ra) \in \HTPQ.$$
That is, $\HTP{R_{(\P-A_0)}}$ is $1$-reducible to $\HTPQ$ for all semilocal
$R_{(\P-A_0)}$, uniformly in $A_0$.
\qed\end{prop}
The proof in \cite{EMPS15}, using work from \cite{K15},
actually shows how to compute, for every prime $p$,
a polynomial $g_{p}(Z,X_1,X_2,X_3)$ such that for all rationals $q$, we have
$$ q\in R_{(\P-\{p\})} \iff g_{p}(q,\Xvec)\in\HTPQ.$$
Therefore, an arbitrary $f(Z_0,\ldots,Z_n)$ has a solution in $R_{(\P-A_0)}$ if and only if
$$(f(\Zvec))^2 + \sum_{p\in A_0, j\leq n} (g_{p}(Z_j,X_{1,j,p},X_{2,j,p},X_{3,j,p}))^2$$
has a solution in $\Q$.
}

\section{Diophantine Undefinability of the Jump}
\label{sec:Kurtz}

The full result of Matiyasevich, Davis, Putnam and Robinson says that not only
is $\HTPZ$ undecidable, but in fact the Halting Problem $\emptyset'$ is \emph{diophantine
in $\Z$}, or expressible in $\Z$ by a diophantine equation.  That is, there exists a polynomial
$f\in\Z[X,Y_1,\ldots,Y_n]$ (for some $n$)
such that
$$ (\forall x\in\Z)~[x\in\emptyset'\iff f(x,Y_1,\ldots,Y_n)\in\HTPZ].$$
Likewise, for any other ring $R$, the sets $S$ \emph{diophantine in $R$} are those
subsets of $R$ definable in the same way by a polynomial in $R[X,Y_1,\ldots,Y_n]$
for some $n$.  When dealing with subrings $R$ of $\Q$, we usually consider only
subsets of $\Z$, often of $\omega$,
using a computable bijection between $\Q$ and $\omega$ if needed.

In our collection of subrings of $\Q$, $\Z$ is the subring $R_\emptyset$, and one
naturally asks whether the proof above carries over to all $W\subseteq\P$:
is $W'$ always diophantine in $R_W$?  Of course, a positive answer would immediately
prove the undecidability of $\HTPQ$, by taking $W=\P$.  In fact, though,
the answer is quickly seen to be negative.  Indeed, with a little more work,
we will show it to be negative in almost all cases.

The easy negative answers arise from taking the set $W$ to be computably
enumerable but not computable.  For example, $W$ might be the image
in $\P$ of $\emptyset'$ itself, under the computable bijection from $\omega$
onto $\P$.  We would then apply the following basic result.  (The reader may
wish to recall the notion of a \emph{$1$-reduction from $A$ to $B$},
which is a computable total function $h:\omega\to\omega$ such that all $x\in\omega$
satisfy $(x\in A\iff h(x)\in B)$.)
\begin{lemma}
\label{lemma:ce}
If $W$ is a computably enumerable subset of $\P$, then $\HTP{R_W}$ is also c.e., and therefore
$\HTP{R_W}\leq_T \emptyset'$.  Indeed there exists a $1$-reduction, i.e.,
$\HTP{R_W}\leq_1 \emptyset'$.

More generally, when $W$ is c.e.\ relative to an oracle $V$, $\HTP{R_W}\leq_1 V'$.
\end{lemma}
\begin{pf}
In the more general setting, let $\la W_s\ra_{s\in\omega}$ be a $V$-computable
enumeration of $W$.  Then $\HTP{R_W}$ is just the set of those $f\in\Z[\Xvec]$
for which
$$(\exists s)(\exists\xvec,y\in\Z)\left[f\left(\frac{x_1}y,\ldots,\frac{x_n}y\right)=0
~\&\text{~all primes dividing $y$ lie in~}W_s\right].$$
Thus $\HTP{R_W}$ is defined by a condition existential relative to $V$ (since
$V$ can compute every $W_s$ uniformly in $s$), and so $\HTP{R_W}$ is
$V$-c.e.  It is then standard (see \cite[Thm\ III.2.3(iii)]{S87})
that $\HTP{R_W}\leq_1 V'$.
\qed\end{pf}

So in the case where $W\subseteq\P$ is the image of $\emptyset'$
(under our computable bijection from $\omega$ onto $\P$), we have
$\HTP{R_W}\leq_T\emptyset'$ but $W'=\emptyset''\not\leq\emptyset'$,
and therefore $W'\not\leq_T \HTP{R_W}$.  Indeed, this holds whenever
the set $W$ is c.e.\ but \emph{non-low}.  Recall that a set $W$ is \emph{low}
if $W'\leq_T\emptyset'$:  this holds of all computable sets $W$, of course,
but also of certain noncomputable sets $W$, both c.e.\ and otherwise.
However, many non-low c.e.\ sets exist, including $\emptyset'$ itself,
and all of these satisfy $W'\not\leq_T\HTP{R_W}$.

Generalizing this result requires theorems proven by Jockusch in \cite{J81}
and by his student Kurtz in his Ph.D.\ thesis \cite{K81}.
\begin{thm}[Jockusch \& Kurtz]
\label{thm:JK}
The set
$$\set{S\subseteq\omega}{(\exists V<_T S)~S\text{~is c.e.\ relative to~}V}$$
has measure $1$ in $2^{\omega}$ (Kurtz) and is comeager there, in the
sense of Baire category (Jockusch).
\end{thm}
A set $S$ in this collection is said to be \emph{relatively c.e.}, as it is c.e.\ in $V$
but not computable in $V$.  Every set $S$ is c.e.\ relative to itself, of course,
but the condition that $S\not\leq_T V$ implies that the jumps of these sets
satisfy $S'\not\leq_1 V'$.  (This follows from the Jump Theorem; see, e.g.,
\cite[Theorem III.2.3]{S87}.)

We can now apply Lemma \ref{lemma:ce} to each relatively c.e.\ set $S$,
with image $W\subseteq\P$, to see that
$$W'\not\leq_1 \HTP{R_W}.$$
Indeed, with $V$ as in the definition, the lemma proves $\HTP{R_W}\leq_1 V'$,
whereas $W'\not\leq_1 V'$ as noted above.  In particular, this yields our next corollary.
\begin{cor}
\label{cor:jump}
The set of those $W\subseteq\P$ such that the jump $W'$ is \emph{not}
diophantine in $R_W$ has full Lebesgue measure within $2^{\P}=\SubQ$, the space
of all subrings of $\Q$, and is comeager there in the sense of Baire category.
\end{cor}
\begin{pf}
If $W'$ is diophantine in $R_W$, via an $f\in\Z[X,\Yvec]$ such that
$$(\forall x\in\omega)~[x\in W'\iff f(x,Y_1,\ldots,Y_n)\in\HTP{R_W}],$$
then we have a $1$-reduction from $W'$ to $\HTP{R_W}$,
by mapping each $x\in\omega$ to the polynomial $f(x,Y_1,\ldots,Y_n)$.
But the discussion above shows that such $1$-reductions exist only for
a measure-$0$, meager class of sets $W\subseteq\P$.
\qed\end{pf}

So the situation for $\Z$ is an anomaly among the subrings of $\Q$.
This is not too surprising:  $\Z$ is very far from being a generic subring
of $\Q$, in any sense of the word ``generic.''  Nevertheless, it is good
to understand that most subrings of $\Q$ do not behave the same way as $\Z$.

It is natural to ask whether one could extend the above result.  In the original
example, with $W$ as the image of $\emptyset'$ in $\P$, we had not only
$W'\not\leq_1\HTP{R_W}$, but actually $W'\not\leq_T\HTP{R_W}$:
there was no Turing reduction, let alone a $1$-reduction.
This would follow more generally for those sets $W$ which are c.e.\
relative to some oracle $V$ such that $W'\not\leq_T V'$.  (One might call such
a $W$ \emph{relatively non-low-c.e.})  However, these sets
are far less common:  the class of subsets of $\omega$ which are
relatively non-low-c.e.\ in this sense has measure $0$ in Cantor space,
and is meager there.  Of course, this does not automatically mean
that $\HTP{R_W}$ must compute $W'$ for the remaining sets $W$ either.
(If it did, then Corollary 1 from \cite{MCiE16} would yield the dramatic
conclusion that $\emptyset'\leq_T\HTPQ$.)

\section{Number Theory}
\label{sec:qpolys}

Our principal tool for proving Theorem \ref{thm:main} and its corollaries,
the chief remaining results in this article, will be the equations
$X^2+qY^2=1$.  In this section we prove the relevant number-theoretic results.
First we show that for each odd prime $q$, there is an
infinite decidable set $V$ of primes such that $R_V$ contains no nontrivial solutions
to $X^2+qY^2=1$.  (Here the trivial solutions are $(\pm 1,0)$, which in Section \ref{sec:thm}
will be ruled out as solutions, at the cost of turning $(X^2+qY^2-1)$ into a messier polynomial.)

\begin{defn}
\label{defn:q-appropriate}
For a fixed odd prime $q$, a prime $p$ is \emph{$q$-appropriate}
if $p$ is odd and $p\neq q$ and $\Leg{-q}p=1$ (that is, $-q$ is a square modulo $p$).
\end{defn}

The crux of Lemmas \ref{lemma:qpolys} and \ref{lemma:p} is that the $q$-appropriate primes
are precisely the possible factors of the denominator in a nontrivial solution
to $x^2+qy^2=1$, thus justifying the term \emph{$q$-appropriate}.

\begin{lemma}
\label{lemma:qpolys}
Fix an odd prime $q$, and let $x$ and $y$ be positive rational numbers with $x^2+qy^2=1$.
Then every odd prime factor $p$ of the least common denominator $c$ of $x$ and $y$ must be $q$-appropriate.

For $q\equiv 3\bmod 4$,
$p$ is $q$-appropriate if and only if $p$ is a square modulo $q$.

For $q\equiv 1\bmod 4$, the situation is a little more complicated.
Now a prime $p$ is $q$-appropriate if and only if one of the following holds:
\begin{itemize}
\item
$p\equiv 1\bmod 4$ and $p$ is a square modulo $q$.
\item
$p\equiv 3\bmod 4$ and $p$ is not a square modulo $q$.
\end{itemize}
\end{lemma}
\begin{pf}
Suppose that 
$a$, $b$, $c$ are positive integers, with no common factor, satisfying $a^2+qb^2=c^2$.
If $p$ divides $c$, then it cannot divide $b$ (lest it also divide $a$),
and so $\left(\frac{a}{b}\right)^2\equiv -q\bmod p$.
Thus every such $p$ is $q$-appropriate.  Suppose in addition that $q\equiv 3\bmod 4$. 
If $p\equiv 1\bmod 4$, then $-1$ is also a square mod $p$, so $q$ is a square mod $p$,
and by quadratic reciprocity $p$ must be a square mod $q$.
On the other hand, if $p\equiv 3\bmod 4$, then $-1$ is not a square mod $p$,
so $q$ is not either; but with both $p$ and $q$ congruent to $3\bmod 4$,
quadratic reciprocity now shows that $p$ is again a square mod $q$.
(The number-theoretic results here may be found in any standard text on the subject,
e.g., \cite{S73}.)

When $q\equiv 1\bmod 4$, a similar analysis, with careful use
of quadratic reciprocity, gives the result stated in the lemma.
\qed\end{pf}

\begin{cor}
\label{cor:q_e seq}
Let $3=q_0<q_1<\cdots$ be the odd prime numbers.
Then, for every $e\in\omega$, there are infinitely
many primes $p$ that are $q_e$-appropriate but (for all $i<e$)
are not $q_i$-appropriate.
\end{cor}
\begin{pf}
A famous theorem of Dirichlet (see \cite[Chap.\ 6, \S 4]{S73}) states that every
arithmetic progression $\set{m+kn}{k\in\omega}$ with $m$ and $n$
relatively prime contains infinitely many primes.
Therefore, the corollary holds for $e=0$, as Lemma \ref{lemma:qpolys}
shows that all primes congruent to $1\bmod 3$ are $3$-appropriate.
Lemma \ref{lemma:qpolys} also noted that the situation is more complicated
for $q_1=5$ than for $q_0=3$, because $q_1\equiv 1\bmod 4$.  Therefore the
full property required for our inductive hypothesis (below) is that all primes
$p\equiv 5\bmod 12$ are $3$-inappropriate. 

Now assume inductively that there is some residue $n$ modulo the product $(4q_0\cdots q_{e-1})$
for which $n\equiv 1\bmod 4$ and no $q_i$ with $i<e$ divides $n$, and such that every prime $p$
with residue $n \bmod (4q_0\cdots q_{e-1})$ is $q_i$-inappropriate for all $i<e$.
(For $e=1$, we saw above that $n=5$ works.)
There are $q_{e}$ distinct elements $n+k(4q_0\cdots q_{e-1})$ in $\Z/(4q_0\cdots q_{e})$,
and their residues modulo $q_{e}$ are all distinct, hence include
all of the elements of $\Z/(q_{e})$.  Lemma \ref{lemma:qpolys}
shows that there are $\frac{q_{e}-1}2$
residues $m$ in $\Z/(4q_0\cdots q_{e})$ such that all
primes with residue $m$ there will be $q_{e}$-appropriate
but $q_i$-inappropriate for all $i<e$.  By Dirichlet's theorem,
for each such $m$, the arithmetic progression $\set{ m+k(4q_0\cdots q_{e})}{k\in\omega}$
will contain infinitely many primes, proving the Corollary for $e$.
On the other hand, another $\frac{q_{e}-1}2$ distinct
residues $m$ in $Z/(4q_0\cdots q_{e})$ have the property
that all primes with that residue are $q_i$-inappropriate for all $i<e+1$,
and that none of $q_0,\ldots,q_{e}$ divides $m$.  Therefore the
inductive hypothesis still holds for $e$, allowing the induction to proceed.

It should be noted that the factor of $4$ in $(4q_0\cdots q_e)$ allowed us to
avoid the bifurcation in Lemma \ref{lemma:qpolys}.  With
the relevant residue $n$ equivalent to $1\bmod 4$, we know that, for all primes
$p$ with that residue, $q_e$-appropriateness simply means being a square modulo $q_e$.
\qed\end{pf}

In fact, the solutions of $X^2+qY^2=1$ are precisely the pairs of the form
$$ \left( \pm\frac{m^2-qn^2}{m^2+qn^2},~\pm\frac{2mn}{m^2+qn^2}\right)$$
for relatively prime integers $m,n\in\N$, not both zero, with the trivial solutions
$(\pm 1,0)$ corresponding to $m=0$ and to $n=0$.  Up to sign, the rational
nontrivial solution $(\frac{a}c,\frac{b}c)$ arises from the integers $m=a-c$ and $n=b$.
For each prime $q$, we call a solution $(a,b,c)$ to $X^2 + qY^2 = Z^2$ {\em primitive}
if $a,b,c$ are pairwise relatively prime positive integers.

\begin{lemma}
\label{lemma:p}
Suppose that $p$ and $q$ are odd primes and $p$ is $q$-appropriate.
Then there is a primitive solution $(a,b,p^k)$ to $X^2 + qY^2 = Z^2$
with $k \ge 1$.  Hence there is a nontrivial solution to $X^2+qY^2 = 1$
in the ring $\Z[\frac1p]$.
\end{lemma}
\begin{pf}
Since $-q$ is a square mod $p$, the ideal $(p)$ in the ring of integers $\cO$ of $K = \Q[\sqrt{-q}]$ splits into distinct prime factors, $(p) = \gp \ov{\gp}$.  If $d$ is the order of $\gp$ in the ideal class group of $\cO$, then $\gp^d = (\alpha)$ is a principal ideal.  The norm $N\!: \, K \to \Q$, given by $N(x+y\sqrt{-q}) = x^2+qy^2$ for $x,y \in \Q$, is multiplicative and so $N(\alpha) = p^d$.   In case $q \equiv 1 \bmod{4}$, the ring of integers $\cO = \Z[\sqrt{-d}]$. Then $\alpha = a+b\sqrt{-q}$ with $a,b \in \Z$ and we are done. 

If $q \equiv 3 \bmod{4}$, then $\cO = \Z[\frac{1 + \sqrt{-q}}{2}]$.  Suppose that $\alpha = \frac{a+b\sqrt{-q}}{2}$ for odd integers $a, b$.  We have $a^2+b^2q = 4p^d$, so  $1 + q \equiv 4 \bmod{8}$ and thus $q \equiv 3 \bmod{8}$.  An elementary computation then shows that $\alpha^3 \in \Z[\sqrt{-q}]$, so we replace $\alpha$ by $\alpha^3$.   In either case, note that $\alpha$ is not divisible by $p$ in $\cO$, so we obtain a primitive solution to $X^2 + qY^2 = Z^2$ of the required form.

To get fancier, one can use the Chebotarev density theorem to choose $p$ so that it splits completely in the Hilbert class field of $\Q[\sqrt{-q}]$.  Then the ideal $\gp = (\alpha)$ already is principal and we can take $d = 1$ above.  (However, if $q \equiv 3 \bmod{8}$, we might still need $\alpha^3$ to get rid of the $2$ in the denominator.)
We refer the reader to \S5 and \S9 of \cite{C} for more information.
\qed\end{pf} 

\comment{
\begin{pf}
Since $-q$ is a square mod $p$, the ideal $(p)$ in the ring of integers $\cO$ of $\Q[\sqrt{-q}]$ splits into distinct prime factors, $(p) = \gp \ov{\gp}$.  If $d$ is the order of $\gp$ in the ideal class group of $\cO$, then $\gp^d = (\alpha)$ is a principal ideal.  By taking the norm $N\!: \, \cO \to \Z$, we obtain $N(\alpha) = p^d$.   In case $q \equiv 1 \bmod{4}$, the ring of integers $\cO = \Z[\sqrt{-d}]$ and we are done. 

If $q \equiv 3 \bmod{4}$, then $\cO = \Z[\frac{1 + \sqrt{-q}}{2}]$.  Suppose that $\alpha = \frac{a+b\sqrt{-q}}{2}$ with $a,b$ odd.  We have $a^2+b^2q = 4p^d$, so  $1 + q \equiv 4 \bmod{8}$ and thus $q \equiv 3 \bmod{8}$.  An elementary computation then shows that $\alpha^3 \in \Z[\sqrt{-q}]$, so we replace $\alpha$ by $\alpha^3$.   Note that $\alpha$ is not divisible by $p$ in $\cO$, so $\alpha$ leads to a primitive solution to $X^2 + qY^2 = Z^2$ of the required form.

To get fancier, one can use the Chebotarev density theorem to choose $p$ so that it  splits completely in the Hilbert class field of $\Q[\sqrt{-q}]$.  Then the ideal $\gp = (\alpha)$ already is principal and we can take $d = 1$ above.  (However, if $q \equiv 3 \bmod{8}$ we might still need $\alpha^3$ to get rid of the 2 in the denominator.)
\qed\end{pf} 
}

\section{Turing Inequivalence}
\label{sec:thm}

Now we apply Lemmas \ref{lemma:qpolys} and \ref{lemma:p} to study the
$HTP$ operator.  It is already known that two Turing-inequivalent sets $U$ and $V$
can have $\HTP{R_U}\equiv_T\HTP{R_V}$.  For a standard example, let $U=\emptyset$
and $V=\emptyset'$:  then $\HTP{R_U}$ has Turing degree $\bfz'$, by the result of
Matiyasevich, Davis, Putnam, and Robinson; whereas
$\HTP{R_V}$ is a c.e.\ set (hence $\leq_T\emptyset'$) which computes $V$,
and thus also has degree $\bfz'$.  The situation here, with $V\equiv_T U'$,
is the maximum possible difference between sets $U$ and $V$
with Turing-equivalent HTP's, since one always has
$V\leq_T\HTP{R_V}\equiv_T\HTP{R_U}\leq_T U'$ and vice versa.
The situation with the jump operator is similar in that Turing-inequivalent sets
$U$ and $V$ can have Turing-equivalent jumps, as with the low
noncomputable sets discussed in Section \ref{sec:Kurtz}.  However, the difference
between the sets $U$ and $V$ could not be a full jump.  (That is,
$U'\not\leq_T V$ and $V'\not\leq_T U$ whenever $U'\equiv_T V'$.)

The jump operator does preserve Turing reducibility, as discussed
earlier.  In fact, $U\leq_T V$ if and only if there is a
$1$-reduction from $U'$ to $V'$, by the standard computability result
known as the Jump Theorem (see e.g.\ \cite[Theorem III.2.3]{S87}).  In contrast,
we now prove that the $HTP$ operator does not preserve Turing reducibility.
Indeed, we will construct a set $U$ for which $\HTP{R_{\Ubar}}$
is Turing-equivalent to the jump $\HTP{R_U}'$.  Once again, this is the
maximum possible difference, since $V\leq_T U$ implies
$$\HTP{R_V}\leq_T V'\leq_T U'\leq_T\HTP{R_U}',$$
with the final reduction holding because $U\leq_T\HTP{R_U}$.
The strong equivalence between our set $U$ and its complement $\Ubar$
makes this all the more striking:  the two sets are Turing-equivalent
via a bounded-truth-table reduction of norm $1$.  (One might pursue
this further, asking whether computably isomorphic sets $U$ and $V$
must have $\HTP{R_U}\equiv_T\HTP{R_V}$.  The set $U$ we build will
not be computably isomorphic to its complement, and this question remains open.)

\begin{thm}
\label{thm:main}
There is a computably enumerable subset $U$ of $\P$,
with complement $\Ubar$ in $\P$, for which $\HTP{R_{\Ubar}}$
computes the set $\textbf{Fin}=\set{e}{W_e\text{~is finite}}$, and therefore
$\HTP{R_{U}}'\equiv_T\HTP{R_{\Ubar}}$.
\end{thm}
\begin{pf}
With $U$ c.e., we will have $\HTP{R_U}\leq_T\emptyset'$, and the Jump
Theorem then shows that $\HTP{R_U}'\leq_T\emptyset''$.
The second jump $\emptyset''$ is well known to satisfy
$\emptyset''\equiv_T \Fin$ (see, e.g., \cite[Theorem IV.3.2]{S87}).
Conversely, $\HTP{R_{\Ubar}}\leq\HTP{R_U}'$,
as discussed above.  Therefore, it is sufficient for us to enumerate a set $U$,
effectively, so that $\Fin\leq_T\HTP{R_{\Ubar}}$.

We accomplish this by using the polynomials $X^2+q_eY^2-1$,
where $q_e$ is the $e$-th odd prime.  Of course,
it is desirable to exclude the trivial solutions $(\pm 1,0)$, so in fact we define
$$ f_e(X,Y,\Zvec,\Tvec) = (X^2+q_eY^2-1)^2+(Y(Z_1^2+\cdots+Z_4^2+1)-(T_1^2+\cdots+T_4^2+1))^2.$$
The second square forces $Y$ to be a quotient of positive rationals,
hence positive.  (With a slightly more complicated polynomial, we could
allow negative values of $Y$ as well, but this is unnecessary here.)
Conversely, the Four Squares Theorem shows that every
positive rational $y$ can be expressed as such a quotient with
all $z_i$ and $t_i$ lying in $\Z$.  Therefore, for every subring $R$ of $\Q$,
$f_e$ lies in $\HTP{R}$ just if $R$ contains elements $x$ and $y>0$
for which $x^2+q_ey^2=1$.

As seen in Section \ref{sec:qpolys}, it is decidable which odd primes $p$ can be factors
of the common denominator of an $x$ and a $y$ with $x^2+q_ey^2=1$.
Lemma \ref{lemma:qpolys} showed that every such $p$ has $\Leg{-q_e}p=1$.
(Recall that we defined such a prime $p$ to be \emph{$q_e$-appropriate}.)
Lemma \ref{lemma:p} proved the converse in a strong way, establishing
that whenever $\Leg{-q_e}p=1$, the ring $\Z[\frac1p]$ already contains
such an $x$ and $y$.  We will enumerate $2$ into our set $U$ immediately,
leaving all odd primes $p$ as candidates for $\Ubar$.  Then we will
enumerate the $q_e$-appropriate primes $p$ into $U$ (that is, out of $\Ubar$)
one by one, as we discover new elements of the c.e.\ set $W_e$.
The goal is that, if $W_e$ is infinite, then all $q_e$-appropriate primes
should be removed from $\Ubar$, so that $\HTP{R_{\Ubar}}$ will not contain $f_e$.
(This goal will not quite be fully achieved, but we will come close enough
to make our proof work.)  However, at any given stage $s$, some particular
$q_e$-appropriate prime $p_{e,s}$ will be \emph{protected by $q_e$},
meaning that for the sake of $q_e$, we will keep $p_{e,s}$ in $\Ubar$
as long as $W_e$ does not acquire any more elements.  If $W_e$ is indeed
finite, then some particular $p_{e,s}$ will be protected from some stage $s_0$
onwards, and therefore will lie in $\Ubar$, forcing $f_e$ to lie in $\HTP{R_{\Ubar}}$.
This will allow our decision procedure for $\Fin$ below an oracle for $\HTP{R_{\Ubar}}$
to succeed.

We use the standard computable enumeration
$\langle W_{e,s}\rangle_{e,s\in\omega}$ of all computably enumerable sets.
This enumeration has the property that for each $s$, there is exactly one
$e$ with $W_{e,s+1}\neq W_{e,s}$, and that for this one $e$, $W_{e,s+1}$
contains all elements of $W_{e,s}$ and exactly one more element as well.

At stage $0$, we set $U_0=\{ 2\}$, and (for every $e\geq 0$) define $q_e$ to be the
$e$-th odd prime and $p_{e,0}$ to be the least $q_e$-appropriate prime
which does not lie in $\{ p_{0,0},\ldots,p_{e-1,0}\}$.
At each stage $s$, the prime $p_{e,s}$ is said to be \emph{protected} by $q_e$ at stage $s$,
although the choice of a protected prime
may change from one stage to the next.  This $p_{e,s}$ is the
prime which $q_e$ currently desires to keep in $\Ubar$.  

At stage $s+1$, we find the unique $e$ with $W_{e,s+1}\neq W_{e,s}$.
This stage is evidence that this particular $W_e$ may be infinite, and so we
enumerate $p_{e,s}$ into $U_{s+1}$.  
For each $i<e$, we keep $p_{i,s+1}=p_{i,s}$.  For each $j\geq e$
(in increasing order), we choose $p_{j,s+1}$ to be the
least $q_j$-appropriate prime which is not in $U_{s+1}$ and
which does not lie in $\{ p_{0,s+1},\ldots,p_{j-1,s+1}\}$.
This clearly makes $p_{e,s+1}\neq p_{e,s}$, since $p_{e,s}\in U_{s+1}$,
but it may leave many subsequent $p_{j,s+1}$ equal to $p_{j,s}$:
the only reason why $p_{j,s+1}$ might not equal $p_{j,s}$ (for $j>e$)
is if $p_{j,s}$ has now been chosen as $p_{k,s+1}$ for some $e\leq k<j$,
i.e., if $p_{j,s}$ is now protected by a higher-priority $q_k$.  This completes
stage $s+1$, and we define the c.e.\ set $U=\cup_s U_s$.

We can describe the arc of a single odd prime $p$ through this process.
Certain primes $p$ might never be protected by any $q_e$:  such a $p$ will lie
in $\Ubar$.  If at some stage $p$ becomes protected, say $p=p_{e,s}$,
then three things can happen.  If no higher-priority $q_i$ (that is, with $i<e$)
subsequently decides to protect $p$, then either
$W_e$ eventually receives a new element and
enumerates $p$ into $U$, or else $W_e$ never receives any new elements
and $p$ stays in $\Ubar$.  The third possibility is that
some higher-priority $q_i$ does subsequently protect $p$, at a stage $s'>s$,
so that $p=p_{i,s'}$ but now $p\neq p_{e,s'}$.  In this case, the same analysis
now applies with $q_i$ in place of $q_e$.  This $p$ could subsequently
be protected by yet another $q_{i'}$ with $i'<i$, but the protecting index can only
change finitely often, of course, since it decreases every time.  Thus, if $p$ ever becomes protected
by any $q_e$, it will either wind up in $U$ or else become the limiting value $p_i$
for some $i\leq e$.

We now prove, by induction on $e$, that the sequence $\la p_{e,s}\ra_{s\in\omega}$
stabilizes on a limit $p_e\notin U$ if $W_e$ is a finite set, but increases without bound
if $W_e$ is infinite.  Let $F=\set{i<e}{W_i\text{~is finite}}$.  First suppose that $W_e$ is finite.
By induction, we may fix a stage $s_0$ such that for all
$i\in F$, the value of $p_{i,s_0}$ never changes at
stages $s>s_0$, and we may also assume $s_0$ to be sufficiently large
that $W_{e,s_0}=W_e$.  By Corollary \ref{cor:q_e seq}, there exist infinitely
many primes which are $q_e$-appropriate but are $q_i$-inappropriate
for all $i<e$.  Such primes will never be
enumerated into $U$ by any $q_i$ with $i<e$ at any stage, and only finitely many
of them can lie in the finite set $U_{s_0}$.  Since the prime $p_{e,s}$
is always chosen to be the least available $q_e$-appropriate prime not already
protected by a higher-priority $q_i$, one of these infinitely many primes will
eventually be chosen as $p_{e,s}$ (unless the sequence $\la p_{e,s}\ra_{s\in\omega}$ stabilizes
on some other prime), and from that stage on, $q_e$ will continue to protect
that same prime $p_{e,s}$:  no higher-priority requirement enumerates it into $U$
because it is $q_i$-inappropriate for all $i<e$;
$q_e$ itself will not enumerate it into $U$ because $W_e$ never again receives a new element;
and once this $p_{e,s}$ has been selected, no $q_j$ with $j>e$ ever again chooses
it as $p_{j,s}$, so no lower-priority $q_j$ ever enumerates it into $U$.
Therefore the limiting value $p_e$ exists as required and lies in $\Ubar$.

On the other hand, suppose $W_e$ is infinite, and fix a stage $s_0$ after which, for all
$i \in F$, the value of $p_{i,s_0}$ never changes again.
Then at each of the infinitely many subsequent stages at which $W_{e,s+1}\neq W_{e,s}$,
$q_e$ will enumerate the current $p_{e,s}$ into $U_{s+1}$ and will choose a new $p_{e,s+1}$.
This $p_{e,s+1}$ is always the least $q_e$-appropriate prime not yet
in $U$ and not currently protected
by any $q_i$ with $i<e$.  Now each $i\in F$ has $p_{i,s}=p_{i,s_0}$
at all these stages, and this limiting value $p_i$ will lie in $\Ubar$.  Some of these
finitely many primes may be $q_e$-appropriate, and so there may be as many as $e$
$q_e$-appropriate primes in $\Ubar$.  However, apart from these $p_i$, every
$q_e$-appropriate prime will eventually be enumerated into $U$.  (If not, then
the least $q_e$-appropriate prime not among these protected $p_i$'s and not in $U$
can never have been protected by any $q_i$ with $i<e$, since such an $i$ would not lie in $F$,
and the least $i$ such that $q_i$ ever protected $p$ would therefore have eventually put
that $p$ into $U$.  But since $p$ was not protected by any $i<e$, it will have
been chosen as $p_{e,s}$ once all smaller $q_e$-appropriate primes not protected
by higher-priority requirements have been enumerated into $U$, and therefore it
too will have been enumerated into $U$.)
Therefore, $p_{e,s}$ does increase without bound as $s\to\infty$, as claimed.
Moreover, with at most $e$ exceptions, all $q_e$-appropriate primes lie in $U$.

Now we give a procedure which uses an $\HTP{R_{\Ubar}}$-oracle to compute $\Fin$.
The procedure decides, for each $e=0,1,2,\ldots$ in turn, whether $W_e$ is infinite or not.
For $e=0$, this is easy, since there is no higher-priority $q_i$ than $q_0$.
If $W_0$ is infinite, $U$ contains all $q_0$-appropriate primes,
and so $f_0\notin\HTP{R_{\Ubar}}$.  On the other hand, if $W_0$ is finite,
then $U$ contains only finitely many $q_0$-appropriate primes, and in particular
does not contain $p_0=\lim_sp_{0,s}$.  By Lemma \ref{lemma:p}, the subring
$\Z[\frac1{p_0}]$ of $R_{\Ubar}$ contains a solution to $f_e$, and so $f_0\in\HTP{R_{\Ubar}}$.

The procedure now continues by recursion on $e$, having determined the finite set $F$ of values $i<e$
lying in $\Fin$.  When it reaches $e$,
it runs the enumeration of $U$
until it finds a stage $s$ for which every $p_{i,s}$ with $i\in F$ lies in $\Ubar$.
(Recall that $\Ubar\leq_T\HTP{R_{\Ubar}}$, so our oracle allows the procedure
to determine membership in $\Ubar$.)  Such a stage must exist, since every such $W_i$
is finite.  For the least $i_0\in F$, $p_{i_0,s}$ can never have become $p_{i',s'}$ for any
$i'<i_0$ at any $s'>s$, since then it would have entered $U$; thus $p_{i_0,s}=p_{i_0}$.
But then the same argument shows inductively for each $i\in I$ that $p_{i,s}=p_i$.
Therefore the procedure uses Proposition \ref{prop:semilocal}  to find a polynomial
$g_e$ which lies in $\HTP{R_{\Ubar}}$ if and only if $f_e\in\HTP{R_{\Ubar-\set{p_{i,s}}{i\in F}}}$.
The oracle then reveals whether $g_e\in\HTP{R_{\Ubar}}$, which in turn determines
whether $e\in\Fin$, since $f_e$ has a solution in the subring $R_{\Ubar-\set{p_i}{i\in F}}$
just if $W_e$ was finite.

It now follows that $\Fin\leq_T\HTP{R_{\Ubar}}$, and therefore $\HTP{R_{\Ubar}}$
has Turing degree at least $\bfz''$, the degree of $\Fin$.  In fact, $\HTP{R_{\Ubar}}$
has exactly this degree, since it must be computable from $\HTP{R_U}'$,
and with $U$ computably enumerable, $\HTP{R_U}$ is c.e.\ as well, forcing
$\HTP{R_U}\leq_T\emptyset'$ and thus $\HTP{R_U}'\leq_T\emptyset''$.
\qed\end{pf}

In the proof of Theorem \ref{thm:main}, our only concern was
to make $\HTP{R_{\Ubar}}$ compute $\emptyset''$, while keeping $U$ c.e.
Now we consider ways to augment this construction.
First, it is not be difficult to use a further finite-injury procedure
to enhance the construction so that it satisfies requirements
to ensure that the c.e.\ set $U$ be \emph{HTP-generic}, in addition to satisfying
$\HTP{R_U}\leq_T\emptyset'$ and $\emptyset''\leq\HTP{R_{\Ubar}}$.
The concept of HTP-genericity is defined and fully explained in \cite[Defn.\ 2]{MCiE16}.
Roughly, it means that for every polynomial $f$, there is a finite initial segment
of $U$ which either ensures that $f\in\HTP{R_U}$ or else ensures that $f\notin\HTP{R_U}$.
It then follows that $\HTP{R_U}\equiv_T U\oplus\HTPQ$.  However, it would not be possible
to make the complement $\Ubar$ HTP-generic:  each $e\notin\Fin$ yields a polynomial
$f\notin\HTP{R_{\Ubar}}$ that has solutions in infinitely many rings $\Z[\frac1p]$,
and therefore no finite subset of the complement of $\Ubar$ (i.e., of $U$) suffices
to guarantee that $f\notin\HTP{R_{\Ubar}}$.  (In the construction above, the $f$
in question is the polynomial $g_e$ for this $e$, since it is necessary
to rule out the finitely many primes $p_i$ with $i<e$ and $i\in\Fin$.)
Indeed, $\Ubar$ must be HTP-nongeneric in order to satisfy $\HTP{R_{\Ubar}}>_T\emptyset'$
and $\Ubar\leq_T\emptyset'$.

Computability theorists familiar with \emph{high permitting} will see a further
enhancement to Theorem \ref{thm:main}: one can make $U$ Turing-reducible
to any given high c.e.\ set $C$.  This is more delicate, and we explain it in
Corollary \ref{cor:main} below, although the proof will be intelligible mainly to those
already familiar with permitting arguments.  Finally, it is not difficult to make $C\leq_1 U$, by a
coding argument that requires only a straightforward finite-injury process.
Therefore, the fully decorated version of Theorem \ref{thm:main} is as follows.

\begin{cor}
\label{cor:main}
For every c.e.\ set $C$ of high degree (i.e., with $\emptyset''\leq_T C'$),
there exists a c.e.\ set $U\subseteq\P$ with $U\equiv_T C$ such that
$$\HTP{R_U}\equiv_T U\oplus\HTPQ\leq_T \emptyset'~~~~\&~~~~\HTP{R_{\Ubar}}\equiv_T \emptyset''.$$
\end{cor}
\begin{pf}
Without rewriting the entire construction, we give reasonable details
about the new condition of high permitting, referring the reader to
\cite[Lemma 12.7.5]{C04} for background.  Given any high c.e.\ degree $\bfc$, we can use
high permitting below $\bfc$ to guarantee that the set $U$ constructed by the theorem
satisfies $\deg{U}\leq_T \bfc$.  Indeed, there must be a c.e.\ set
$C\in\bfc$ with a computable enumeration $\la C_s\ra_{s\in\omega}$
such that the computation function of $C$ using this enumeration dominates every
total computable function.  High permitting
requires that we allow a prime $p$ to enter $U$ only at a stage $s+1$
such that $C_{s+1}\res p\neq C_s\res p$.  So, when $W_{e,s+1}\neq W_{e,s}$
we may not be allowed to enumerate the current $p_{e,s}$ into $C$ immediately.
Instead, we \emph{mark} the prime for $U$, meaning that we
put it on a waiting list to enter $U$.  Whenever a number $m$ appears
in $C_{t+1}-C_t$, all primes $> m$ currently on the waiting list are enumerated
into $U_{t+1}$.  This guarantees that $U\leq_T C$, since $p\in U$ just if $p\in U_s$,
where $s$ is the least stage for which $C_s\res p=C\res p$ (and this stage $s$
can be computed given a $C$-oracle).

The principal difference from Theorem \ref{thm:main} is that with high permitting,
in the case where $e\in\Inf$, we do not know exactly how many $q_e$-appropriate primes
will have to be left out of $U$, although the total number left out will be finite.
Therefore, instead of knowing exactly which question to ask about $\HTP{R_{\Ubar}}$
to determine whether an $e$ lies in $\Inf$, we will need to resort to a search for a
finite set $A\subseteq\P$ such that $f_e\notin\HTP{R_{\Ubar-A}}$, employing Corollary
\ref{cor:semilocal}.  This will yield an enumeration of $\Inf$, and we will then apply
Lemma \ref{lemma:Inf} (below).  In turn, we will need to ensure, whenever $e\in\Fin$,
not just that one particular limit prime $p_e$ remains in $\Ubar$, but that infinitely
many $q_e$-appropriate primes lie in $\Ubar$; otherwise we would mistakenly enumerate
$e$ into $\Inf$.  The full requirement is:
$$\mathcal{R}_e:  e\in\Inf\iff\text{$\Ubar$ contains only finitely many $q_e$-appropriate primes.}$$
Therefore, for each $e$ at stage $s$, we do not define just one protected prime $p_{e,s}$,
but rather make a list $p_{e,0,s}<p_{e,1,s}<\cdots$ of all $q_e$-appropriate primes
not yet marked for $U$ nor protected by any higher-priority requirement at that stage.
The basic rule of protection is that $\R_{e+1}$ cannot mark $p_{e,0,s}$ for $U$ (but can mark
any other $p_{e,j,s}$), $\R_{e+2}$ cannot mark $p_{e,0,s}$ nor $p_{e,1,s}$,
and in general $\R_{e+k}$ cannot mark any $p_{e,j,s}$ with $j<k$ for $U$.
By the same token, $\R_e$ itself cannot mark $p_{e-1,0,s}$, nor $p_{e-2,0,s}$ nor $p_{e-2,1,s}$, etc.,
so these primes will not be chosen as $p_{e,i,s}$.  This protection rule ensures that $\R_e$ will
not force any lower-priority $\R_i$ to leave infinitely many $q_i$-appropriate primes in $\Ubar$,
but also that $\R_e$ will have a choice of infinitely many $q_e$-appropriate primes to keep in $\Ubar$
if necessary, since at some stage $s$ some $q_{e+1}$-inappropriate prime will be chosen
as $p_{e,1,s}$, and later some $p_{e,2,s}$ will be chosen which is both $q_{e+1}$-inappropriate
and $q_{e+2}$-inappropriate, and so on.  Indeed, if $e\in\Fin$, then every limit
$p_{e,i}=\lim_s p_{e,i,s}$ will exist and all these limits $p_{e,i}$ will be
$q_e$-appropriate and will lie in $\Ubar$.  On the other hand, if $e\in\Inf$,
then all but finitely many
$q_e$-appropriate primes will be chosen as $p_{e,i,s}$ at some stage
(or else enumerated into $U$ by a higher-priority requirement), and thus
all but finitely many $q_e$-appropriate primes will eventually
be marked for $U$, since the current $p_{e,0,s}$ is marked for $U$ whenever $W_e$
gets a new element.

Next, it is necessary to see that for each $e\in\Inf$, only finitely many
of the $q_e$-appropriate primes ever marked for $U$ fail to enter $U$.
To see this, notice that if $p$ is marked for $U$ at a stage $s$ but never
enters $U$, then $C_s\res p=C\res p$.  That is, $s\geq C_C(p)$, where
$C_C$ is the computation function defined (as in \cite[p.\ 230]{C04})
so that $C_C(x)$ is the least $s$ with $C_s\res x=C\res x$.
By hypothesis, this $C_C$ is not only noncomputable, but dominates
every total computable function.  If $e\in\Inf$, then as seen above,
cofinitely many $q_e$-appropriate primes are eventually marked for $U$ by $q_e$.
Let these primes be $p_0<p_1<\cdots$, which is a computable infinite sequence,
and define the computable function $f(n)$ to equal the stage at which
$p_n$ is marked for $U$.  Since $C_C$ dominates $f$, there are only
finitely many $n$ with $f(n)\geq C_C(n)$, and for all other $n$ than these,
we have $C_{f(n)}\res n\neq C\res n$, by the definition of $C_C$.
But $p_n>n$, so also $C_{f(n)}\res p_n\neq C\res p_n$, and therefore,
at some stage $s$ after the stage $f(n)$ at which $p_n$ is marked for $U$,
$C$ will permit $p_n$ to enter $U$.  Therefore $\R_e$ is indeed satisfied.

This being the case, Proposition \ref{prop:semilocal} now allows us
to use an oracle for $\HTP{R_{\Ubar}}$ to enumerate the set
$\Inf$, the complement of $\Fin$.  To do so, for each $e\in\omega$,
we simply go through the finite initial segments $A_n=\{ 2,3,5,\ldots,p_n\}$
of $\P$, one set at a time.  (One could speed up this process by
considering only initial segments of the set of $q_e$-appropriate primes.)
For each $n$, we ask the $\HTP{R_{\Ubar}}$-oracle whether the polynomial
$$ (f_e(X,Y))^2 + \sum_{i\leq n} (g_{p_i}(X,Z_1,Z_2,Z_3))^2 + \sum_{i\leq n} (g_{p_i}(Y,T_1,T_2,T_3))^2$$
has a solution in $R_{\Ubar}$.  If the answer is ever negative, then we know that
$e\in\Inf$, because all $q_e$-appropriate primes except those in that $A_n$
must lie in $U$.  On the other hand, if no $n$ ever yields a negative answer, then
infinitely many $q_e$-appropriate primes must lie in $\Ubar$, and so $e\in\Fin$.
Therefore, Lemma \ref{lemma:Inf} below proves that $\Inf\leq_T\HTP{R_{\Ubar}}$.
Adding the coding of $C$ into $U$ is standard, as is the addition of requirements
for HTP-genericity, since both of these are finitary and blend easily with the permitting.
\qed\end{pf}

The next result follows directly from Corollary \ref{cor:main}, but is nevertheless
quite striking:  the $HTP$ operator can reverse strict Turing reducibility.

\begin{cor}
\label{cor:leapfrog}
There exist subrings $R$ and $S$ of $\Q$ with
$R<_T S$, yet with $\HTP{S}<_T\HTP{R}$.
\end{cor}
\begin{pf}
Using an incomplete high c.e.\ set $C$ in Corollary \ref{cor:main} yields
a set $\Ubar<\emptyset'$, yet $\HTP{R_{\emptyset'}}\equiv_T\emptyset'$, being c.e.,
while $\emptyset''\leq\HTP{R_{\Ubar}}$.  So we let $R=R_{\Ubar}$ and $S=R_{\emptyset'}$.
(Alternatively, as the c.e.\ sets are dense under Turing reducibility, we could
fix a c.e.\ set $D$  with $C<_T D<_T\emptyset'$ and apply
Corollary \ref{cor:main} to $D$, getting $S=R_V$ for a c.e.\ set $V\equiv_T D$, and
coding $\emptyset'$ into $\HTP{R_V}$ as well, thus making $S<_T\HTP{S}$.)
\qed\end{pf}

It remains to establish the lemma required for Corollary \ref{cor:main},
which is a standard computability result.
\begin{lemma}[Folklore]
\label{lemma:Inf}
For subsets $A\subseteq\omega$,
$$ \Inf \leq_T A \iff \Inf\text{~is $A$-computably enumerable.}$$
\end{lemma}
\begin{pf}
The forward direction is immediate.  For the converse, we define a computable
total function $f$ so that, for every $e$,
$$\phi_{f(e)}(s) =\left\{\begin{array}{cl} 0,&\text{~if~}\phi_{e,s}(e)\diverges;\\
\diverges,&\text{~if~}\phi_{e,s}(e)\converges.
\end{array}\right.$$
Thus $e$ lies in the complement
$\overline{\emptyset'}$ of the Halting Problem if and only if the domain $W_{f(e)}$ of $\phi_{f(e)}$
is infinite.  (That is, $f$ is an \emph{$m$-reduction} from $\overline{\emptyset'}$ to $\Inf$.)
But now, since $\Inf$ is $A$-computably enumerable, so is $\overline{\emptyset'}$,
and therefore, with an $A$-oracle, we can compute
$\emptyset'$.  Since $\Fin$ is c.e.\ relative to $\emptyset'$, and $\Inf$
is already $A$-c.e., this allows us to compute $\Inf$ using an $A$-oracle.
\qed\end{pf}

High permitting, as opposed to ordinary c.e.\ permitting, does appear necessary
in Corollary \ref{cor:main}.  With ordinary permitting, one could ensure
that infinitely many $q_e$-appropriate primes entered $U$, in the case where
$e\in\Inf$, but the requirement that cofinitely many should enter $U$ requires
the strength of high permitting.  If we could have enumerated an HTP-generic
$U$ below a non-high $C$, while still achieving $\emptyset''\leq_T\HTP{R_{\Ubar}}$,
then $\HTP{R_U}$ would have been high, hence $>_T U$, which (with $U$
HTP-generic) would have established the undecidability of $\HTPQ$.
The fact that this is not possible here essentially means that, while $\HTPQ$
may yet turn out to be undecidable, the polynomials $X^2+qY^2=1$
are not complex enough to prove it.

As a final side note, it is possible to adjust the construction in Theorem \ref{thm:main}
to ensure that $\emptyset''\leq_1\HTP{R_{\Ubar}}$.  Since $\HTP{R_{\Ubar}}\leq_1 U'$,
this situation is only possible when $\emptyset'\leq_T U$, and thus cannot be accomplished
using high permitting below an incomplete set $C$.  The adjustment simply stipulates
that, for all $e<i$, if $p$ is the least prime which is $q_j$-appropriate for all $e\leq j\leq i$,
then the only $q_i$-appropriate primes which $q_e$ ever restrains from
entering $U$ should be those which are $\leq p$.  When $q_e$ redefines its protected prime $p_{e,s+1}$,
it takes this stipulation into account, and immediately enumerates into $U$ all
$q_e$-appropriate primes in the interval $[p_{e,s},p_{e,s+1})$ which are
not protected by higher-priority requirements.

\section{Questions}
\label{sec:questions}

There is a natural analogy between the $HTP$ operator, mapping $W$ to $\HTP{R_W}$,
and the \emph{jump operator}, mapping $W$ to $W'$.  $W'$ and $\HTP{R_W}$ are both
$W$-computably enumerable, and as noted earlier, the basic situation for Turing reducibility is that
$$ W\leq_T \HTP{R_W} \leq_T W',$$
with equality possible at either end, though of course not at both ends simultaneously.

The analogy is strengthened by the parallels between $\HTPQ$ and the Halting
Problem $\emptyset'$.  The class $\GL1$ of \emph{generalized low$_1$ sets} $W$
is defined by the property
$$ W' \equiv_T \emptyset'\oplus W.$$
This class is comeager of measure $1$ in Cantor space.  Of course,
$\emptyset'\oplus W\leq_1 W'$ always holds.  The opposite reduction is trickier:
it does fail on a meager set of measure $0$, and even within $\GL1$ it is
in general only a Turing reduction, not a $1$-reduction.  This opposite reduction
holds uniformly on a comeager class, but not on any class of measure $1$.
That is, there is a single Turing functional $\Phi$ such that
$$\set{W\subseteq\omega}{\chi_{W'}=\Phi^{\emptyset'\oplus W}}$$
is comeager, but for every $\Phi$, this class fails to have measure $1$.  (Given $\ep>0$,
one can choose a $\Phi$ for which it has measure $>1-\ep$, and the choice
of the program for $\Phi$ is uniform in $\ep\in\Q$.)

For the $HTP$ operator, the analogy leads us to the $HTP$-generic sets,
introduced at the end of the previous section.  All $HTP$-generic sets $W\subseteq\P$
satisfy
$$ \HTP{R_W} \equiv_T \HTPQ\oplus W.$$
The class of HTP-generic sets is comeager in Cantor space, but its measure there is unknown.
Again, $\HTPQ\oplus W\leq_1\HTP{R_W}$ always holds, while in the opposite direction,
Turing reducibility holds uniformly on a comeager class.  (Details about these results
appear in \cite{MCiE16}.)

On the other hand, as the name implies, the $HTP$-generic sets are defined
by a genericity property.  In this sense, they are analogous to the \hbox{\emph{$1$-generic sets}:}
those $U\subseteq\omega$ such that, for every $e$,
$$\mu(\set{V\in 2^\omega}{(\exists n)~V\res n \in W_e}=1\implies (\exists n)~U\res n\in W_e.$$
The $1$-generic sets are precisely the points
in $2^\omega$ at which the jump operator is continuous, and likewise the HTP-generic
sets are precisely the points of continuity of the HTP operator.  All $1$-generic sets lie
in $\GL1$, and the $1$-generic sets form a comeager class in $2^\omega$, but of measure $0$.
Our analogy therefore suggests that the class of HTP-generic sets, already
known to be comeager, might also have measure $0$.

As one can see, the state of knowledge about $HTP$-genericity is strong \textit{vis-\`a-vis}
Baire category, but less so \textit{vis-\`a-vis} Lebesgue measure.  In work yet to appear,
the author has shown that if the $HTP$-generic sets have measure $1$, then there can be no
existential definition of $\Z$ in the field $\Q$.  Theorem \ref{thm:main} here is an existence
theorem, but says nothing about measure; indeed, the construction used in its proof
involves (for certain values $e$, namely those in $\Inf$) enumerating cofinitely many
$q_e$-appropriate primes into $U$.  Clearly, even for a single $q$, the class
of sets $U\subseteq\P$ which contain all but finitely many of the $q$-appropriate primes is
a class of measure $0$, and so Theorem \ref{thm:main} yields no conclusions
about measure.  New results regarding measure would be of real interest for the study
of definability and decidability in number theory.

\parbox{4.7in}{
{\sc
\noindent
Department of Mathematics \hfill \\
\hspace*{.1in}  Queens College -- C.U.N.Y. \hfill \\
\hspace*{.2in}  65-30 Kissena Blvd. \hfill \\
\hspace*{.3in}  Queens, New York  11367 U.S.A. \hfill \\
\ \\
Ph.D. Program in Mathematics \hfill \\
Ph.D. Program in Computer Science (Miller) \hfill \\
\hspace*{.1in}  C.U.N.Y.\ Graduate Center\hfill \\
\hspace*{.2in}  365 Fifth Avenue \hfill \\
\medskip
\hspace*{.3in}  New York, New York  10016 U.S.A. \hfill}\\
\hspace*{.045in} {\it E-mail: }
\texttt{kkramer\at {qc.cuny.edu} }\hfill \\
\hspace*{.045in} {\it E-mail: }
\texttt{Russell.Miller\at {qc.cuny.edu} }\hfill \\

\comment{
{\sc
\noindent
Department of Mathematics \hfill \\
\hspace*{.1in}  Queens College -- C.U.N.Y. \hfill \\
\hspace*{.2in}  65-30 Kissena Blvd. \hfill \\
\hspace*{.3in}  Queens, New York  11367 U.S.A. \hfill \\
Ph.D. Programs in Mathematics \& Computer Science \hfill \\
\hspace*{.1in}  C.U.N.Y.\ Graduate Center\hfill \\
\hspace*{.2in}  365 Fifth Avenue \hfill \\
\hspace*{.3in}  New York, New York  10016 U.S.A. \hfill}\\
\medskip
\hspace*{.045in} {\it E-mail: }
\texttt{Russell.Miller\at {qc.cuny.edu} }\hfill \\
}
}

\end{document}